%abstract in the end
\input amstex
\loadeufm
\loadeusm
\loadbold

%\input amstex
%\documentstyle{amsppt}
%\NoRunningHeads
%\hsize=36 true pc
%\vsize=54 true pc
%\voffset=0.08 true in
%\magnification = \magstep1

\documentstyle{amsppt}
\NoRunningHeads
\magnification=1200
\pagewidth{32pc}
\pageheight{42pc}
\vcorrection{1.2pc}

\define\wh{\widehat}

\define\bc{\Bbb C}
\define\bz{\Bbb Z}

\define\ep{\epsilon}
\define\vep{\varepsilon}

\define\gln{U_q(\widehat{\frak{gl}}_{{}_N})}

\define\sln{U_q(\widehat{\frak{sl}}_{{}_N})}

\define\oh{\overline{h}}
\define\tor{U_q(\frak{sl}_{{}_N, tor})}

\topmatter
\title  $U_q(\widehat{\frak{gl}}_{{}_N})$ action on $\widehat{\frak{gl}}_{{}_N}$-modules and
  quantum toroidal algebras
\endtitle
\footnotetext"$^{1}$"{The author gratefully acknowledges the grant support from the
Natural Sciences and Engineering Research Council of Canada.}
\footnotetext"$^{2}$"{The author gratefully acknowledges the support
from  NSF Grant DMS-9970493.}

\author Yun Gao
\footnotemark"$^{1}$" 
and Naihuan Jing
\footnotemark"$^{2}$"
\endauthor

\address
Department of Mathematics and Statistics, York University, Toronto, \newline
Canada M3J 1P3
\endaddress
\email ygao\@yorku.ca
\endemail
\address
Department of Mathematics, North Carolina State University, Raleigh,
NC 27695-8205, USA
\endaddress
\email jing\@math.ncsu.edu
\endemail

\endtopmatter
\document

\subhead\S 0. Introduction\endsubhead

\medskip

A quantum Kac-Moody algebra $U_q(\frak{g})$ introduced in [D1] and [Ji] is a
$q$-deformation of the universal enveloping algebra $U(\frak{g})$ of a Kac-Moody Lie
algebra $\frak g$. Lusztig [L] (see also [R]) has shown that the character formula for
dominant highest
weight representations is preserved when $U_q(\frak{g})$ is deformed to $U(\frak{g})$.

This paper deals with the quantum affine algebra $U_q(\widehat{\frak{gl}}_{{}_N})$,
the quantum toroidal algebra $\tor$, and vertex operators. Our purpose is to give an irreducible 
vertex representation for the newly developed quantum toroidal
algebra $U_q(\frak{sl}_{{}_N, tor})$. This leads an interesting phenomena: the irreducible basic $\widehat{\frak{gl}}_{{}_N}$-module allows an irreducible
$U_q(\widehat{\frak{gl}}_{{}_N})$ action. Moreover, each weight space as both 
$U_q(\widehat{\frak{gl}}_{{}_N})$-module and $\widehat{\frak{gl}}_{{}_N}$-module coincides. 
Therefore,
this fact enhances the above mentioned Lusztig theorem in the special case
$\frak{g}=\widehat{\frak{gl}}_{{}_N}$.

Quantum toroidal algebras were  introduced by Ginzburg-Kapranov-Vasserot
  [GKV] in the study of the Langlands reciprocity for algebraic
surfaces. These
algebras
are quantized analogues for toroidal Lie algebras of Moody-Rao-Yokonuma [MRY]. 
A Schur type duality between representations of
the quantum toroidal algebra (of type $A$) and 
the double affine Hecke algebra 
was established by Varagnolo-Vasserot in [VV1]. 
They further obtained 
in [VV2](see also [STU]) a nice  representation 
for the quantum toroidal algebra by gluing the standard 
module of the quantum 
affine algebra together with a level-$0$ module arising from solvable lattice models.
Some other interesting representations were constructed in [FJW], [S] and [TU] from 
various point of view. 

Vertex representations for quantum affine algebras were developed by Frenkel-Jing [FJ] and
Jing [J], which are q-analogues of Frenkel-Kac [FK] and Segal [Se] construction for affine
Lie algebras. Vertex representations of the affine Lie algebra $\widehat{\frak gl}_N$
was given in Frenkel \cite{F}. Representations for $U_q(\widehat{\frak{gl}}_{{}_N})$ or
$U_q(\widehat{\frak{sl}}_{{}_N})$ have been studied by many people (for example, [H], [JKKMP] and 
[KMS]). 

We shall use  $N$ copies of orthogonal (or independent)
Fock spaces to 
construct
a  family of vertex operators as was done in [G]. More precisely, given a non-zero complex number $p$, 
 we define vertex operators 
$X_{ij}(r, z)$ depending on the parameter $p$, for $i, j, r\in\Bbb Z$. 
In the non-quantum case (i.e. $q=1$), there are only finite such
operators for a fixed $r$ as $X_{ij}(r, z)$ is doubly periodic with respect to indices $i$ and $j$. However,
in the quantum case, due to the $q$-twisting, $X_{ij}(r, z)$ is no longer doubly periodic 
so there are infinite such operators
for a fixed $r$. The algebra generated by some of those operators will give an irreducible representation
of the quantum toroidal algebra $U_q(\frak{sl}_{{}_N, tor})$. As one by-product, we obtain another
realization of the quantum affine algebra $U_q(\widehat{\frak{gl}}_{{}_N})$
 and recover the 
 vertex representation of [FJ] for the quantum affine algebra $U_q(\widehat{\frak{sl}}_{{}_N})$.
The other by-product is the fact of level-0 $U_q(\widehat{\frak{gl}}_{{}_N})$ action on
level-1 $U_q(\widehat{\frak{gl}}_{{}_N})$-module which appeared in [JJKMP].

Our results are based on a new realization of the quantum affine algebra
$\gln$. In \cite{DF} the quantum affine algebra $\gln$ was defined in the 
context of the quantum Yang-Baxter equation, which contains the quantum affine
algebra $\sln$ in a canonical way. 
We construct a new set of generators within
the quantum affine algebra $\gln$ with the advantage that they are orthogonal
(mutually commutative) in the usual sense, and this immediately
 establishes the isomorphism between
the classical enveloping algebra of $\widehat{\frak gl}_N$ and its quantum
counterpart over the field of functions in $q^{c/2}$. Drinfeld \cite{D2}
showed that
in an appropriate completion the quantum enveloping algebra of the
simple Lie algebra is isomorphic to the enveloping algebra of the simple
Lie algebra up to a twisting of the Hopf algebra structure. 
Our result adds another example to this general
statement and also facilitates the construction of Weyl bases in the algebra. 

The main results in sections 1 to 3 were announced in [GJ]. The paper is organized as follows. In Sect. 1 we recall the notion of the quantum toroidal algebra. Sect. 2 lays
the foundation for the Fock space representation of
$\tor$ where we emphasize the role of the affine general linear algebra. Sect. 3 is devoted to the proof of our construction. In the last section we construct a new set of
generators for the quantum affine algebra $\gln$ and through the new basis we show that the module provides both
actions for $\gln$ and $\tor$.

\medskip

\subhead \S 1. Quantum toroidal algebra of type $A_{{}_{N-1}}$ \endsubhead

\medskip

%The definition of the quantum toroidal algebra is given by using the %formal power series which
%is a two dimensional analog of Drinfeld realization of quantum affine %algebras.
%Recall that $\delta(z)= \sum_{k\in \Bbb Z} z^k.$

We always assume that the complex number $q$ is generic
and $N$ is a positive integer with $N\geq 3$. Let $d$ be a nonzero complex number. 
The quantum 
toroidal algebra $U_q(\frak{sl}_{{}_N, tor})$ is the unital associative algebra over $\Bbb C$
generated by $e_{i, k}, f_{i, k}, h_{i, l}, k_i^{\pm 1}$, where $i=0, 1, \cdots, N-1, 
k\in \Bbb{Z}, l\in \Bbb{Z}\setminus\{ 0\}$, and the central elements $c^{\pm 1}$. The relations
are expressed in terms of the formal series:
$$ e_i(z) = \sum_{k\in \Bbb{Z}}e_{i, k}z^{-k}, \quad f_i(z) = \sum_{k\in \Bbb{Z}}f_{i, k}
z^{-k}, \tag 1.1$$
and 
$$k_i^{\pm}(z) = k_i^{\pm 1} \exp(\pm(q- q^{-1})\sum_{k=1}^{\infty} h_{i, \pm k} z^{\mp k}),\tag 1.2 $$
as follows
$$\align & k_i k_i^{-1} = k_i^{-1}k_i = c c^{-1}=1, \tag 1.3\\
&[k_i^{\pm}(z), k^{\pm}_j(w)] = 0,\tag 1.4 \\
&\theta_{-a_{ij}}(c^2 d^{-m_{ij}}\frac{w}{z})k_i^+(z)k_j^-(w) =
 \theta_{-a_{ij}}(c^{-2}d^{-m_{ij}}\frac{w}{z})k_j^-(w)k_i^+(z), \tag 1.5\\
&k_i^{\pm}(z) e_j(w) = 
\theta_{\mp a_{ij}}(c^{-1}d^{\mp m_{ij}}(\frac{w}{z})^{\pm 1})e_j(w)k_i^{\pm}(z), \tag 1.6\\
&k_i^{\pm}(z) f_j(w) = 
\theta_{\pm a_{ij}}(c d^{\mp m_{ij}}(\frac{w}{z})^{\pm 1})f_j(w)k_i^{\pm}(z), \tag 1.7\\
&[e_i(z), f_j(w)] =\frac{\delta_{ij}}{q-q^{-1}}(\delta(c^{-2}\frac{z}{w}) k_i^+(cw)
-\delta(c^2\frac{z}{w})k_i^-(cz)), \tag 1.8\\
&(d^{m_{ij}}z - q^{a_{ij}}w)e_i(z) e_j(w) = (q^{a_{ij}}d^{m_{ij}}z - w)e_j(w) e_i(z) , \tag 1.9
\endalign$$
$$\align
&(q^{a_{ij}}d^{m_{ij}}z - w)f_i(z) f_j(w) = (d^{m_{ij}}z - q^{a_{ij}}w)f_j(w) f_i(z) , \tag 1.10\\
&\{e_i(z_1) e_i(z_2)e_j(w) - (q+q^{-1})e_i(z_1)e_j(w)e_i(z_2) + e_j(w)e_i(z_1)e_i(z_2)\}
 \tag 1.11\\
& + \{ z_1 \leftrightarrow z_2\} = 0, \text{ if } a_{ij} = -1, 
\endalign$$
$$
\align
&\{f_i(z_1) f_i(z_2)f_j(w) - (q+q^{-1})f_i(z_1)f_j(w)f_i(z_2) + f_j(w)f_i(z_1)f_i(z_2)\}
 \tag 1.12 \\
& + \{ z_1 \leftrightarrow z_2\} = 0, \text{ if } a_{ij} = -1, \\
&[e_i(z), e_j(w)] = [f_i(z), f_j(w)]= 0,  \text{ if } a_{ij} = 0, \tag 1.13
\endalign$$
where 
$$\theta_m(z)=\frac{q^m z-1}{z-q^m} \in \Bbb{C}[[z]] \tag 1.14 $$ is understood as
the Taylor series expansion, 
$$ A=(a_{ij}) =\pmatrix 2 & -1 &  & 0 & -1  \\
-1 & 2 & \ldots & 0 & 0 \\
&  \vdots & \ddots & \vdots & \\
0 & 0 & \ldots & 2 & -1 \\
-1 & 0 & & -1 & 2 \endpmatrix \tag 1.15 $$ 
is the Cartan matrix of
affine type 
$A^{(1)}_{{}_{N-1}}$  
 and 
$$M=(m_{ij})= \pmatrix 0 & -1 &  & 0 & 1  \\
1 & 0 & \ldots & 0 & 0 \\
&  \vdots & \ddots & \vdots & \\
0 & 0 & \ldots & 0 & -1 \\
-1 & 0 & & 1 & 0 \endpmatrix \tag 1.16 $$ 
is a skew symmetric matrix (i.e., $m_{ij}=\delta_{i, j+1}- \delta_{j, i+1}$ for $0\leq i, j\leq N-1$,
where $\delta_{0N}=\delta_{N0}=1$).

The reason to call $U_q(\frak{sl}_{{}_N, tor})$ the quantum toroidal algebra  is that it  is a two-parameter
deformation of the enveloping algebra of the toroidal Lie algebra 
$\hat{\frak{sl}}_{{}_N}(\Bbb{C}[x^{\pm 1}, y^{\pm 1}])$
(the universal central extension of the double loop algebra 
$\frak{sl}_{{}_N}(\Bbb{C}[x^{\pm 1}, y^{\pm 1}])$, 
see [MRY]).
Actually,  $U_q(\frak{sl}_{{}_N, tor})$ has been proven to be  a one-parameter deformation of the enveloping 
algebra of a Lie algebra
over a quantum torus, see Section 13 in [VV2]. 

As pointed out in [GKV], the quantum toroidal algebra $U_q(\frak{sl}_{{}_N, tor})$ contains two remarkable subalgebras,
 the horizontal subalgebra $\dot{U}_h$ and the vertical subalgebra $\dot{U}_v$, where $\dot{U}_h$ is generated by
 $$e_{i, 0}, f_{i, 0}, k_i^{\pm 1} , \quad i = 0, 1, 2, \cdots, N-1, \tag 1.17$$  while 
 $\dot{U}_v$ is generated by 
 $$e_{i, n}, f_{i, n}, h_{i, l}, k_i^{\pm 1}, \quad  i =1, 2, \cdots, N-1, n\in\Bbb{Z} \text{ and } l\in \Bbb{Z}^\times.
 \tag 1.18$$
The central elements of $\dot{U}_h$ and $\dot{U}_v$ are 
$k^{\pm 1}=\prod_{i=0}^{N-1}k_i^{\pm 1}$ and $c^{\pm 1}$ respectively. 
  Both 
 $\dot{U}_h$ and $\dot{U}_v$ are isomorphic to the quantum affine algebra $\sln$.

\medskip

\subhead \S 2. Fock space and vertex operators \endsubhead

\medskip

In this section, we shall set up our Fock space and construct a family of vertex operators indexed
by $\Bbb{Z}\times \Bbb{Z}$.

Let
$$P=\bz\ep_1 \oplus \cdots \oplus\bz\ep_{{}_N}\tag 2.1$$
be a rank $N$ free abelian group provided with a $\bz$-bilinear form $(\cdot ,\cdot)$ defined by $(\ep_i, \ep_j) = \delta_{ij}$, $1\leq i, j\leq N$. 
Let
$$Q=\bz(\ep_1-\ep_2)\oplus \cdots \oplus\bz(\ep_{{}_{N-1}}-\ep_{{}_N})\tag 2.2$$
be the rank $(N-1)$ free subgroup of $P$. Then
$$\Delta =\{\alpha\in Q:  (\alpha, \alpha) = 2\}=
\{\ep_i -\ep_j : 1\leq i \neq j\leq N\} \tag 2.3$$
is the root system of type $A_{{}_{N-1}}$.

Let $\vep: Q\times Q \to \{\pm 1\}$ be a bimultiplicative function 
such that 
$$\align
\vep(\alpha + \beta, \gamma) &= \vep(\alpha, \gamma)\vep(\beta, \gamma), 
\ \ 
\ \vep(\alpha, \beta +\gamma) = \vep(\alpha, \beta)\vep(\alpha, \gamma) \tag 2.4 \\ 
\vep(\alpha, \alpha) &= (-1)^{\frac{1}{2}(\alpha, \alpha)}, \, \, \tag 2.5
\endalign $$
for $\alpha, \beta, \gamma\in Q$.
 The formula 
(2.5) immediately implies
$$\vep(\alpha, \beta)\vep(\beta, \alpha) = (-1) ^{(\alpha, \beta)}, 
\qquad \alpha, \beta \in Q. \tag 2.6$$

Let $\bc[Q] =\sum \oplus\bc e^\alpha $
be the group algebra of $Q$. For $\beta\in Q$, define $e_\beta\in \text{End}\bc[Q]$ by
$$e_\beta e^\alpha = \vep(\beta, \alpha)e^{\alpha +\beta}, \, \text{ for  }\alpha 
\in Q.\tag 2.7$$
It follows  that
$$e_\alpha e_\beta = (-1)^{(\alpha, \beta)}e_{\beta}e_\alpha \tag 2.8$$
for $\alpha, \beta\in Q$. Also, for $\beta\in H = Q\otimes_{\Bbb Z}\Bbb C$, define $\beta(0)\in \text{End}\bc[Q]$ by
$$\beta(0) e^\alpha = (\beta, \alpha) e^\alpha, \, \text{ for } \alpha\in Q.\tag 2.9$$

Next let $\ep_i(n)$ and $C$ be the generators of the Heisenberg algebra $\Cal H$, $ 1\leq i \leq N, n\in\bz\setminus\{0\}$, subject to relations that $C$ is central and
$$
[\ep_i(m), \ep_j(n)]=m\delta_{ij}\delta_{m+n, 0}C.
 \tag 2.10
$$ 

Define
$$\ep_{ij}(n)=q^{(j-i)|n|/2}\ep_{i}(n)-q^{(i-j)|n|/2}\ep_j(n).
 \tag 2.11 $$
Here we observe the relation
$$\ep_{ik}(n)=q^{(k-j)|n|/2}\ep_{ij}(n)+q^{(i-j)|n|/2}\ep_{jk}(n). \tag 2.12$$

We have
$$
[\ep_{ij}(m), \ep_{ij}(n)]=m\delta_{m+n, 0}(q^{(j-i)m}+
q^{(i-j)m})C, \qquad j>i, \tag 2.13
$$
and it is zero if $i=j$.

Let
$$S(\Cal{H}^-) = \bc[\ep_i(n): 1\leq i\leq N, n\in -\bz_+]\tag 2.14$$
denote the symmetric algebra of $\Cal{H}^-$, which is the algebra of polynomials in 
infinitely many variables $\ep_i(n), 1\leq i\leq N, n\in -\bz_+$, where $\bz_+ =\{ n\in \bz: n>0\}$.
$S(\Cal{H}^-)$ is
 an $\Cal H$-module in which $C=1$, $\ep_i(n)$ acts as the 
multiplication operator for $n\in -\bz_+$, and $\ep_i(n)$ acts as the partial
 differential operator for $n\in \bz_+$.

For $\alpha \in \{\pm\ep_1, \cdots, \pm\ep_N\}$
we introduce the operators $E_{\pm}(r, z)$ as follows
$$
E_{\pm}(\alpha, z)=
\exp(\mp\sum_{n=1}^{\infty}\frac{\alpha(\pm n)}{n}z^{\mp n}). 
\tag 2.15
$$
It follows that for  $\alpha, \beta\in \{\pm\ep_1, \cdots, \pm\ep_N\}$
$$
E_+(\alpha, z)E_-(\beta, w) = E_-(\beta, w)E_+(\alpha, z) 
(1-\frac{w}{z})^{(\alpha, \beta)}. \tag 2.16
$$

Set 
$$V_Q = S(\Cal{H}^-)\otimes \bc[Q].\tag 2.17$$

The operator $z^\alpha \in (\text{End}\bc[Q])[z, z^{-1}]$ is defined as 
$z^{\alpha(0)}=\exp(\alpha(0)\ln z)$.
$$z^\alpha e^\beta = z^{(\alpha, \beta)}e^\beta\tag 2.18$$
for $\alpha, \beta\in Q$. Thus, in $(\text{End}\bc[Q])\{z\}$, we have
$$[\alpha(0), z^\beta] = 0, \quad \text{ and } \quad z^\alpha e_\beta = 
e_\beta z^{\alpha + (\alpha, \beta)} \tag 2.19$$
for $\alpha, \beta\in Q$. It is clear that the formula (2.18) expresses 
$z^\alpha$ as an operator from $\bc[Q]$ to $\bc[Q][z, z^{-1}]$.

Let $\mu$ be any non-zero complex number. Consider the valuation $\mu^\alpha$ of 
the operator $z^\alpha$. Namely, $\mu^\alpha$ is the operator $\bc[Q]\to \bc[Q]$ given by
$$\mu^\alpha e^\beta = \mu^{(\alpha, \beta)} e^\beta, \, \text{ for } \alpha, \beta\in Q.\tag 2.20$$
%This operator $\mu^\alpha$ has been used in [J] and [G]. It will be useful %in our construction as well.

Most notations related to formal series and the Fock space used in this section 
can be found in [FLM] and [J].

\medskip

Set 
$$\ep_{i+N} = \ep_i, \quad \text{ for } i\in \Bbb Z.$$
Accordingly,
$$(\ep_i, \ep_j) = \delta_{ij} = \delta_{\bar{i}, \bar{j}}, \text{ for } \bar{i}, \bar{j}\in \Bbb{Z}/N\Bbb{Z}.$$
This implies that $a_{ij}=(\ep_i -\ep_{i+1}, \ep_j-\ep_{j+1})$ and $m_{ij}=(\ep_i, \ep_{j+1})
-(\ep_j, \ep_{i+1})$ which are used in Section 1.

Let $p$ ba a non-zero complex number. For $r,  i, j\in \Bbb Z$, we define the vertex operator 
$X_{ij}(r, z)$ as follows.
$$
\align
X_{ij}(r, z) 
=&
:\exp(-\sum_{n\neq 0}\frac{(\ep_{i}(n)-p^{-rn}q^{(i-j)|n|}\ep_{j}(n))}{n}z^{-n}):  \\
&\qquad e_{\ep_i-\ep_j}z^{\ep_i-\ep_j+\frac{(\ep_i-\ep_j, \ep_i -\ep_j)}{2}}p^{-r\ep_j- \frac{(\ep_j, \ep_i -\ep_j)}{2}r} \tag 2.21\\
=& E_-(\ep_i, z)
E_-(-\ep_j, zq^{i-j}p^r)
 E_+(\ep_i, z)E_+(-\ep_j, zq^{j-i}p^r)\\
&\qquad e_{\ep_i-\ep_j}z^{\ep_i-\ep_j+
\frac{(\ep_i-\ep_j, \ep_i -\ep_j)}{2}}p^{-r\ep_j- \frac{(\ep_j, \ep_i -\ep_j)}{2}r}
\endalign
$$

\noindent{\bf Remark 2.22.} In the non-quantum case (i.e. $q=1$), we get
$X_{ij}(r, z) = X_{i+N, j}(r, z)= X_{i, j+N}(r, z).$
 For general $q$,   we only  have
$X_{ij}(r, z) = X_{i+N, j+N}(r, z).$ 
We may think of $X_{10}(r, z)$ (or $X_{01}(r, z)$) as $X_{1N}(r, z)$ (or $X_{N1}(r, z)$)
by modulo $N$ in the powers of $q$.

\medskip

Let
$$\tilde{\frak{gl}}_{{}_N} 
=\frak{gl}_{{}_N}(\Bbb{C}[t, t^{-1}])\oplus\Bbb{C}C\oplus\Bbb{C}D$$
be the affine algebra, where $C$ is the central element and $D$ is the degree operator. 
$\frak{h}=H\oplus \Bbb{C}C\oplus \Bbb{C}D$ is its Cartan subalgebra.

The following result is known (see [F] and [G]).

\proclaim{Proposition 2.23} $V_Q$ is an irreducible $\tilde{\frak{gl}}_{{}_N}$ module given by
$$\align &E_{ij}(t^n) \mapsto X_{ij}(0,n), \quad 1\leq i\neq j\leq N, n\in\Bbb{Z},\\ 
&E_{ii}(t^n)\mapsto \ep_i(n), \quad 1\leq i\leq N, n\in\Bbb{Z}\\
&C\mapsto 1, \quad D\mapsto -\frac{1}{2}\sum_{i=1}^N\ep_i(0)^2  -\sum_{i=1}^N\sum_{n\in\bz_+}
\ep_i(-n)\ep_i(n).
\endalign$$
and
$V_Q= \sum_{\mu\in P}\oplus V_\mu$, 
where 
$P=\{\omega_0 + \alpha + n\tau : \alpha\in Q, 
\, n\in \bz_{\leq 0}\}$
is the set of weights, $V_\mu =\{ v\in V: h.v = \mu(h)v \text{ for }h\in \frak{h}\}$. 
More precisely, we have
$V_{\omega_0 +\alpha + n\tau_0} = 
  W\otimes e^\alpha,$
where $W$ is the homogeneous subspace of $S(\Cal{H}^-)$ of degree 
$n +\frac{1}{2}(\alpha, \alpha)$, for $\alpha\in Q, 
n\in\bz_{\leq 0}$, and  
$$\text{ch } V_Q =(\sum_{\alpha\in Q}e^{\omega_0 + 
\alpha -\frac{1}{2}(\alpha, \alpha)\tau})\varphi(e^{-\tau})^{-N},\tag 2.24$$
where
$\omega_0|_{H\oplus \bc D}=0, \, \, 
\omega_0(C) =1, $
$\tau|_{H\oplus \bc C} = 0,  \, \,  \tau(D) =1.$
\endproclaim

Next, for $r, i, j\in \Bbb Z$, and $i \neq j $,  we define
$$\align & u_{ij}(r, z)
 = q^{(j-i)(\ep_i -\ep_j)} \tag 2.25 \\
 &\cdot \exp(\sum_{n\geq 1}\frac{q^{(j-i)n}-q^{(i-j)n}}{n}
( q^{\frac{j-i}{2}n}\ep_i(n)-p^{-nr}q^{\frac{i-j}{2}n}\ep_j(n)z^{-n})\\
& v_{ij}(r, z) 
= q^{(i-j)(\ep_i -\ep_j)} \tag 2.26 \\
& \cdot \exp(
\sum_{n\geq 1}\frac{q^{(i-j)n}-q^{(j-i)n}}{n}(q^{\frac{j-i}{2}n}\ep_i(-n)-p^{nr}q^{\frac{i-j}{2}n}
\ep_j(-n))z^{n}).
\endalign $$

The normal ordering can be defined as usual. For instance we have
$$\align & :X_{ij}(r_1, z_1)X_{kl}(r_2, z_2): \tag 2.27\\
%=&\exp(-\sum_{n\in -\bz_+}\frac{(\ep_i(n)-p^{-nr_1}q^{(i-%j)|n|}\ep_j(n))z_1^{-n} + 
%(\ep_k(n) -p^{-nr_2}q^{(k-l)|n|}\ep_l(n))z_2^{-n}}{[n]})\\
%&\cdot\exp(-\sum_{n\in \bz_+}\frac{(\ep_i(n)-p^{-nr_1}q^{(i-%j)|n|}\ep_j(n))z_1^{-n} +
% (\ep_k(n) -p^{-nr_2}q^{(k-l)|n|}\ep_l(n))z_2^{-n}}{[n]})\\
%&\cdot e_{\ep_i-\ep_j}e_{\ep_k-\ep_l}\\
%&\cdot z_1^{\ep_i-\ep_j +\frac{(\ep_i-\ep_j, \ep_i-\ep_j +\ep_k-%\ep_l)}{2}}
% z_2^{\ep_k-\ep_l +\frac{(\ep_k-\ep_l, \ep_i-\ep_j +\ep_k-\ep_l)}{2}}\\
%&\cdot p^{-r_1\ep_j -r_2\ep_l -\frac{(\ep_j, \ep_i-\ep_j)}{2}
%r_1-\frac{(\ep_l, \ep_k-\ep_l)}{2}r_2}\\
=&E_-(\ep_i, z_1)E_-(-\ep_j, z_1q^{i-j}p^{r_1})E_-(\ep_k, z_2)E_-(-\ep_l, z_2q^{k-l}p^{r_2})\\
& E_+(\ep_i, z_1)E_+(-\ep_j, z_1q^{j-i}p^{r_1})E_+(\ep_k, z_2)E_+(-\ep_l, z_2q^{l-k}p^{r_2})\\
&\cdot e_{\ep_i-\ep_j}e_{\ep_k-\ep_l}
z_1^{\ep_i-\ep_j +\frac{(\ep_i-\ep_j, \ep_i-\ep_j +\ep_k-\ep_l)}{2}}
 z_2^{\ep_k-\ep_l +\frac{(\ep_k-\ep_l, \ep_i-\ep_j +\ep_k-\ep_l)}{2}}\\
&\cdot p^{-r_1\ep_j -r_2\ep_l -\frac{(\ep_j, \ep_i-\ep_j)}{2}
r_1-\frac{(\ep_l, \ep_k-\ep_l)}{2}r_2}\\
\endalign$$

Moreover, we have 
 $$:X_{ij}(r_1, z_1)X_{kl}(r_2, z_2): = (-1)^{(\ep_i-\ep_j, \ep_k-\ep_l)}
:X_{kl}(r_2, z_2)X_{ij}(r_1, z_1):.\tag 2.28$$
where $r_1, r_2, i, j , k, l \in\Bbb Z$ .

The following fact is straightforward.

\proclaim{Lemma 2.29} For $r_1, r_2,  i,  j, k, l\in \Bbb Z$,
$$\align &X_{ij}(r_1, z_1)X_{kl}(r_2, z_2)=
: X_{ij}(r_1, z_1)X_{kl}(r_2, z_2) : \\
 &\qquad \cdot (\frac{z_1}{z_2})^{\frac{(\ep_i-\ep_j, \ep_k-\ep_l)}{2}}
 (1-\frac{z_2}{z_1})^{\delta_{ik}}
 (1-\frac{p^{r_2}q^{k-l}z_2}{p^{r_1}q^{j-i}z_1})^{\delta_{jl}}   
\tag 2.30
\\
 &\qquad \cdot (1-\frac{p^{r_2}q^{k-l}z_2}{z_1})^{-\delta_{il}}
 (1-\frac{z_2}{p^{r_1}q^{j-i}z_1})^{-\delta_{jk}} p^{-\delta_{jk}r_1 +\delta_{jl}r_1} \endalign$$
\endproclaim

\noindent {\bf Remark 2.31.} Note that $\delta_{ij}$ used above is doubly periodic with respect
to $i$ and $ j$ as $\delta_{i, j}=(\ep_i, \ep_j)= \delta_{\bar{i}, \bar{j}}$, where  $\bar{i}, \bar{j}\in \Bbb{Z}/N\Bbb{Z}$.

To calculate the commutators of vertex operators, we need some more notations 
and identities. Set
$$\align & F_{kl}^{ij}(r_1, r_2, z_1, z_2)=: X_{ij}(r_1, z_1)X_{kl}(r_2, z_2) : \\ 
&\qquad \cdot(\frac{z_1}{z_2})^{\frac{(\ep_i-\ep_j, \ep_k-\ep_l)}{2}}p^{-\delta_{jk}r_1 +\delta_{jl}r_1}   
(1-\frac{z_2}{z_1})^{\delta_{ik}}
 (1-\frac{p^{r_2}q^{k-l}z_2}{p^{r_1}q^{j-i}z_1})^{\delta_{jl}} \tag 2.32\\
 &\qquad \cdot (1-\frac{p^{r_2}q^{k-l}z_2}{z_1})^{1-\delta_{il}}
 (1-\frac{z_2}{p^{r_1}q^{j-i}z_1})^{1-\delta_{jk}} \frac{p^{r_1}z_1}{z_2}.
\endalign$$

In particular we have $F_{ji}^{ij}(r_1, r_2, z_1, z_2) = F_{ij}^{ji}(r_2, r_1, z_2, z_1)=:X_{ij}(r_1, z_1)X_{ji}(r_2, z_2):$. Then we can rewrite
Lemma 2.29 as follows.
$$\align &X_{ij}(r_1, z_1)X_{kl}(r_2, z_2) \\
=& F_{kl}^{ij}(r_1, r_2, z_1, z_2)(1-\frac{p^{r_2}q^{k-l}z_2}{z_1})^{-1}(1-\frac{z_2}
{p^{r_1}q^{j-i}z_1})^{-1}\frac{z_2}{p^{r_1}z_1}. \tag 2.33
\endalign$$

The next lemma can be checked directly.

\proclaim{Lemma 2.34} For $r_1, r_2, i, j\in \Bbb Z$ with $r_1 + r_2 =0 $ and $i \neq j $,
$$\align 
&\lim_{z_1\to p^{r_2}q^{j-i}z_2} F_{ji}^{ij}(r_1, r_2, z_1, z_2)=:X_{ij}(r_1, p^{r_2}q^{j-i}z_2)X_{ji}(r_2, z_2):\\
&\qquad =u_{ij}(r_1, p^{r_2}q^{\frac{j-i}{2}}z_2) \\
&\lim_{z_2\to p^{r_1}q^{j-i}z_1} F_{ji}^{ij}(r_1, r_2, z_1, z_2)=:X_{ij}(r_1, z_1)X_{ji}(r_2, p^{r_1}q^{j-i}z_1):\\
&\qquad =v_{ij}(r_1, q^{\frac{j-i}{2}}z_1) 
\endalign$$
\endproclaim

The following basic result is similar to  (4.16) in [J] whose proof is straightforward.

\proclaim{Lemma 2.35} For any $a, b\in \bc$ and $a\neq b$, we have
in $\bc[[z, z^{-1}]]$
$$(1-az)^{-1}(1-bz)^{-1} = \frac{z^{-1}}{a-b}((1-az)^{-1}-(1-bz)^{-1})$$
\endproclaim

\proclaim{Proposition 2.36} If $r_1 + r_2 = 0$, then as a formal series we have
$$\align &(1-\frac{p^{r_2}q^{j-i}z_2}{z_1})^{-1}(1-\frac{z_2}{p^{r_1}q^{j-i}z_1})^{-1}\frac{z_2}{p^{r_1}z_1}
-(1-\frac{p^{r_1}q^{i-j}z_1}{z_2})^{-1}
(1-\frac{z_1}{p^{r_2}q^{i-j}z_2})^{-1} \frac{z_1}{p^{r_2}z_2}\\
=&(q^{j-i}-q^{i-j})^{-1}(\delta(\frac{p^{r_2}q^{j-i}z_2}{z_1})
-\delta(\frac{p^{r_1}q^{j-i}z_1}{z_2}))
\endalign$$
\endproclaim
\demo{Proof} By  Lemma 2.35, we see that the left hand side of the identity is equal to
$$\align &(q^{j-i} -q^{i-j})^{-1}((1-\frac{p^{r_2}q^{j-i}z_2}{z_1})^{-1}
-(1-\frac{z_2}{p^{r_1}q^{j-i}z_1})^{-1})
\\
&-(q^{i-j}-q^{j-i})^{-1}((1-\frac{p^{r_1}q^{i-j}z_1}{z_2})^{-1}-(1-\frac{z_1}{p^{r_2}q^{i-j}z_2})^{-1})\\
=&(q^{j-i}-q^{i-j})^{-1}(\delta(\frac{q^{j-i}p^{r_2}z_2}{z_1})-
\delta(\frac{q^{j-i}p^{r_1}z_1}{z_2}))
\endalign$$
as needed. \qed\enddemo

\medskip

Now we are in the position to show our commutation relation:

\proclaim{Proposition 2.37} If $r_1+r_2= 0$, then
$$\align &[X_{ij}(r_1, z_1), X_{ji}(r_2, z_2)]  \\
=&\frac{1}{q^{j-i}- q^{i-j}}(u_{ij}(r_1, q^{\frac{j-i}{2}}p^{r_2}z_2)\delta(\frac{q^{j-i}p^{r_2}z_2}{z_1}) 
-v_{ij}(r_1, q^{\frac{j-i}{2}}z_1)\delta(\frac{q^{j-i}p^{r_1}z_1}{z_2}))\\
\endalign$$
\endproclaim
\demo{Proof} By (2.33), we have
$$\align &[X_{ij}(r_1, z_1), X_{ji}(r_2, z_2)] \\
=&X_{ij}(r_1, z_1)X_{ji}(r_2, z_2) - X_{ji}(r_2, z_2)X_{ij}(r_1, z_1)\\
=&F_{ji}^{ij}(r_1, r_2, z_1, z_2)(1-\frac{p^{r_2}q^{j-i}z_2}{z_1})^{-1}
(1-\frac{z_2}{p^{r_1}q^{j-i}z_1})^{-1}\frac{z_2}{p^{r_1}z_1}\\
&-F_{ij}^{ji}(r_2, r_1, z_2, z_1)(1-\frac{p^{r_1}q^{i-j}z_1}{z_2})^{-1}
(1-\frac{z_1}{p^{r_2}q^{i-j}z_2})^{-1}\frac{z_1}{p^{r_2}z_2}.
\endalign$$

From Proposition 2.36 and the normal ordering relations of the F's, the above becomes
$$F_{ji}^{ij}(r_1, r_2, z_1, z_2)(q^{j-i}-q^{i-j})^{-1}(\delta(\frac{q^{j-i}p^{r_2}z_2}{z_1})
- \delta(\frac{q^{j-i}p^{r_1}z_1}{z_2})).$$
Now (2.37) follows from (2.34).  \qed
\enddemo

\medskip

\subhead \S 3. Representations for the quantum toroidal algebra \endsubhead

\medskip

In this section, we shall use some of vertex operators constructed in the previous section
to generate a unital associative
 algebra which turns out to be a homomorphic image of $U_q(\frak{sl}_{{}_N, tor})$.
From now on we  assume that
$$p= d^{-N}. \tag 3.1$$

For simplicity, we let
$$ \align & E_i(z) = X_{i, i+1}(0, p^{\frac{1}{2}}d^iz), \quad 
F_i(z)= X_{i+1, i}(0, p^{\frac{1}{2}}d^i z), \quad i = 1, 2, \cdots, N-1,\\
& E_0(z)= X_{01}(1, p^{-\frac{1}{2}}z), \quad \quad
F_0(z) = X_{10}(-1, p^{\frac{1}{2}}z),\\
& K^+_i(z) = u_{i, i+1}(0, p^{\frac{1}{2}}d^i z), \quad 
K^-_i(z) = v_{i, i+1}(0, p^{\frac{1}{2}}d^i z), \quad i =1, 2, \cdots, N-1,\\
& K_0^+(z) = u_{01}(1, p^{-\frac{1}{2}}z), \quad \quad 
K_0^-(z) = v_{01}(1, p^{-\frac{1}{2}}z),\quad
K^{\pm}_i = q^{\pm(\ep_i -\ep_{i+1})}.
\endalign $$

Let $\Cal A$ be the associative algebra generated by the coefficients of $E_i(z), F_i(z),
K_i^{\pm}(z)$, for $0 \leq i \leq N-1$. 

\medskip

We now state our main result of this paper.

\proclaim{Theorem 3.2} The linear map $\pi$ given by
$$ \align & \pi(e_i(z)) = E_i(z), \quad 
\pi(f_i(z))= F_i(z), \quad i = 1, 2, \cdots, N-1,\\
& \pi(e_0(z))= E_0(z), \quad 
\pi(f_0(z)) = F_0(z),\\
& \pi(k^+_i(z)) = K^+_i(z), \quad 
\pi(k^-_i(z)) = K^-_i(z) , \quad i =1, 2, \cdots, N-1,\\
&\pi(k_0^+(z)) = K_0^+(z), \quad 
\pi(k_0^-(z)) = K_0^-(z) ,\\
& \pi(c) = q^{\frac{1}{2}},  \pi(k_i^{\pm}) = K_i^{\pm}, \quad i=0, 1, \cdots, N-1
\endalign $$
yields  an algebra homomorphism from $U_q(\frak{sl}_{{}_N, tor})$ to $\Cal A$.

\endproclaim

\demo{Proof} The proof will be carried out in several steps.

\medskip

Step 1. First we have
$$u_{i, i+1}(r_1, z)v_{j, j+1}(r_2, w)=
v_{j, j+1}(r_2, w)u_{i, i+1}(r_1, z)\exp(-(q-q^{-1})^2J)
$$
where
$$ J
= \sum_{n\geq 1}
[\frac{q^{\frac{1}{2}n}\ep_i(n)-p^{-nr_1}q^{-\frac{1}{2}n}\ep_{i+1}(n)}{n}, 
\frac{q^{\frac{1}{2}n}\ep_j(-n)-p^{nr_2}q^{-\frac{1}{2}n}\ep_{j+1}(-n)}{n}](\frac{w}{z})^n$$
Hence the above becomes
$$\align 
=& v_{j, j+1}(r_2, w)u_{i, i+1}(r_1, z) \\
&\cdot \frac{(1-\frac{q^3w}{z})^{\delta_{ij}}(1-\frac{q^{-1}w}{z})^{\delta_{ij}}}
{(1-\frac{qw}{z})^{2\delta_{ij}}}
\frac{(1-\frac{p^{r_2-r_1}qw}{z})^{\delta_{ij}}(1-\frac{p^{r_2-r_1}q^{-3}w}{z})^{\delta_{ij}}}
{(1-\frac{p^{r_2-r_1}q^{-1}w}{z})^{2\delta_{ij}}}\\
&\cdot \frac{(1- \frac{p^{r_2}w}{z})^{2\delta_{i, j+1}}}
{(1-\frac{q^{-2}p^{r_2}w}{z})^{\delta_{i,j+1}} 
(1- \frac{q^2p^{r_2}w}{z})^{\delta_{i, j+1}}} 
\frac{(1- \frac{p^{-r_1}w}{z})^{2 \delta_{j, i+1}}}
{(1-\frac{q^{-2}p^{-r_1}w}{z})^{\delta_{j, i+1}} 
(1- \frac{q^2p^{-r_1}w}{z})^{\delta_{j, i+1}}} 
\endalign $$
By restricting $0\leq i, j\leq N-1$ and $r_1, r_2 =0, 1$, noting that $p=d^{-N}$ which is only
involved when $K_0^{\pm}(z)$ presents, we thus obtain
$$\theta_{-a_{ij}}(q d^{-m_{ij}}\frac{w}{z})K_i^+(z)K_j^-(w) =
 \theta_{-a_{ij}}(q^{-1}d^{-m_{ij}}\frac{w}{z})K_j^-(w)K_i^+(z),\tag 3.3$$
 for $0\leq i, j\leq N-1$.

\medskip

Step 2. The relations (1.6-1.7) are proved similarly.
We only show one case in the following.

$$\align &u_{i, i+1}(r_1, z) X_{jk}(r_2,  w)=X_{jk}(r_2, w) u_{i, i+1}(r_1, z)q^{(\ep_i -\ep_{i+1}, \ep_j - \ep_k)}\\
&\cdot \exp((q-q^{-1})\sum_{n\geq 1}[q^{\frac{n}{2}}\ep_i(n)-p^{-nr_1}q^{-\frac{n}{2}}\ep_{i+1}(n),
\frac{\ep_j(-n) - p^{nr_2} q^{(j-k)n}\ep_k(-n)}{[n]}](\frac{w}{z})^n)\\
=&X_{jk}(r_2, w) u_{i, i+1}(r_1, z)q^{(\ep_i -\ep_{i+1}, \ep_j - \ep_k)}\frac{(1-\frac{q^{-\frac{1}{2}}w}{z})^{\delta_{ij}}}
{(1-\frac{q^{\frac{3}{2}}w}{z})^{\delta_{ij}}} 
\frac{(1-\frac{q^{j-k-\frac{3}{2}}p^{r_2-r_1}w}{z})^{\delta_{i+1, k}}}
{(1-\frac{q^{j-k + \frac{1}{2}}p^{r_2 -r_1}w}{z})^{\delta_{i+1, k}}} \\
&\cdot \frac{(1-\frac{q^{j-k + \frac{3}{2}}p^{r_2}w}{z})^{\delta_{ik}}}
{(1-\frac{q^{j-k-\frac{1}{2}}p^{r_2}w}{z})^{\delta_{ik}}}
\frac{(1-\frac{q^{\frac{1}{2}}p^{-r_1}w}{z})^{\delta_{i+1, j}}}
{(1-\frac{q^{-\frac{3}{2}}p^{-r_1}w}{z})^{\delta_{i+1, j}}}
\endalign$$
By restricting $0\leq i , j, k \leq N$, $k=j \pm 1$ and $r_1, r_2 =0, 1$, noting that $p=d^{-N}$ 
which is only
involved when $K_0^{+}(z)$, $E_0(z)$ or $F_0(z)$ presents, we get
for $0\leq i, j\leq N-1$
$$\align &K_i^{+}(z) E_j(w) = 
\theta_{- a_{ij}}(q^{-\frac{1}{2}}d^{- m_{ij}}(\frac{w}{z}))E_j(w)K_i^{+}(z), \tag 3.4\\
&K_i^{+}(z) F_j(w) = 
\theta_{ a_{ij}}(q^{\frac{1}{2}} d^{- m_{ij}}(\frac{w}{z}))F_j(w)K_i^{+}(z). \tag 3.5\\
\endalign$$

Step 3. Assume that $0\leq i \leq j \leq N - 1$. One has
$$\align &X_{i, i+1}(r_1, z) X_{j, j+1}(r_2,  w)=:X_{i, i+1}(r_1, z) X_{j, j+1}(r_2,  w) : \\
& \qquad \cdot (\frac{z}{ w})^{\delta_{ij}-\frac{\delta_{i+1,j}}{2}}(1-\frac{w}{z})^{\delta_{ij}}
(1-\frac{p^{r_2}q^{-1}w}{p^{r_1}qz})^{\delta_{ij}}
(1-\frac{w}{p^{r_1}qz})^{-\delta_{i+1, j}} p^{-\delta_{i+1, j}r_1 + \delta_{ij}r_1}\\
&X_{j, j+1}(r_2, w) X_{i, i+1}(r_1, z)=:X_{j, j+1}(r_2, w) X_{i, i+1}(r_1, z) : \\
&\qquad\cdot (\frac{w}{z})^{\delta_{ij}-\frac{\delta_{i+1, j}}{2}}
(1-\frac{z}{w})^{\delta_{ij}}
(1-\frac{p^{r_1}q^{-1}z}{p^{r_2}qw})^{\delta_{ij}}
(1-\frac{p^{r_1}q^{-1}z}{w})^{-\delta_{i+1, j}} p^{\delta_{ij} r_2}.
\endalign $$

It follows that 
$$\align &
(d^{m_{ij}}z - q^{a_{ij}}w)E_i(z) E_j(w) = (q^{a_{ij}}d^{m_{ij}}z - w)E_j(w) E_i(z) ,\tag 3.6\\
& [E_i(z), E_j(w)] = 0,  \text{ if } a_{ij} = 0. \tag 3.7
\endalign $$
Note that $a_{ij}=0$ if and only if $|i- j|\neq 0, 1$ and $p=d^{-N}$ which is only involved
when $E_0(z)$ presents. We can prove relation (1.10) similarly.

Assume that $0\leq i \neq j \leq N - 1$.
$$\align &X_{i, i+1}(r_1, z) X_{j+1, j}(r_2,  w)=:X_{i, i+1}(r_1, z) X_{j+1, j}(r_2,  w) : \\
& \qquad\cdot (\frac{z}{w})^{\frac{\delta_{i, j+1}+\delta_{i+1,j}}{2}}
(1-\frac{w}{z})^{\delta_{i,j+1}}
(1-\frac{p^{r_2}w}{p^{r_1}z})^{\delta_{i+1,j}}p^{\delta_{i+1, j}r_1}\\
&X_{j+1, j}(r_2, w) X_{i, i+1}(r_1, z)=
:X_{j+1, j}(r_2,w) X_{i, i+1}(r_1, z) : \\
&\qquad\cdot (\frac{w}{z})^{\frac{\delta_{i, j+1} + \delta_{i+1, j}}{2}}
(1-\frac{z}{w})^{\delta_{i, j+1}}
(1-\frac{p^{r_1}z}{p^{r_2}w})^{\delta_{i+1, j}}p^{\delta_{j, i+1}r_2}
\endalign $$

It immediately  gives us
$$[E_i(z), F_j(w)] = 0, \tag 3.8 $$
for $0\leq i\neq j\leq N-1$.
\medskip

Step 4. To prove the Serre relation, we set
$$\align & I(z_1, z_2, w)= \frac{w}{(z_1z_2)^{1/2}}\left(
\frac{(z_1-z_2)(z_1-q^{-2}z_2)}{(z_1-q^{-1}w)(z_2-q^{-1}w)}\right.\\
&\quad + \left. (q+q^{-1})
\frac{(z_1-z_2)(z_1-q^{-2}z_2)}{(z_1-q^{-1}w)(w-q^{-1}z_2)}
+ \frac{(z_1-z_2)(z_1-q^{-2}z_2)}{(w-q^{-1}z_1)(w-q^{-1}z_2)}\right).
\endalign
$$

\proclaim{Lemma 3.9} \cite{J} Let $S(z_1, z_2, w) = I(z_1, z_2, w) + I(z_2, z_1, w)$, then 
$$S(z_1, z_2, w) = 0.$$
\endproclaim
It follows from the usual vertex operator computation that
$$\align & X_{ij}(r_1, z_1)X_{ij}(r_2, z_2)X_{kl}(r_3, w)\\ 
&\quad = :X_{ij}(r_1, z_1)X_{ij}(r_2, z_2)X_{kl}(r_3, w):  \\
&\quad \cdot \frac{z_1}{z_2}(1-\frac{z_1}{z_2})(1-\frac{p^{r_2}z_2}{p^{r_1}z_1})p^{r_1}
(\frac{z_1}{w})^{\frac{(\ep_i-\ep_j, \ep_k-\ep_l)}{2}}
(1-\frac{w}{z_1})^{\delta_{ik}} \\
&\quad \cdot (1-\frac{p^{r_3}q^{k-l}w}{p^{r_1}q^{j-i}z_1})^{\delta_{jl}}
 (1-\frac{p^{r_3}q^{k-l}w}{z_1})^{-\delta_{il}}
 (1-\frac{w}{p^{r_1}q^{j-i}z_1})^{-\delta_{jk}} \tag3.10\\
 &\quad \cdot p^{-\delta_{jk}r_1 +\delta_{jl}r_1}(\frac{z_2}{w})^{\frac{(\ep_i-\ep_j, \ep_k-\ep_l)}{2}}
(1-\frac{w}{z_2})^{\delta_{ik}}(1-\frac{p^{r_3}q^{k-l}w}{p^{r_2}q^{j-i}z_2})^{\delta_{jl}}  \\
 &\quad \cdot
  (1-\frac{p^{r_3}q^{k-l}w}{z_2})^{-\delta_{il}}
 (1-\frac{w}{p^{r_2}q^{j-i}z_2})^{-\delta_{jk}} p^{-\delta_{jk}r_2 +\delta_{jl}r_2}. \\
\endalign$$
Its associated products $X_{ij}(r_1, z_1)X_{kl}(r_3, w)X_{ij}(r_2, z_2)$, $X_{kl}(r_3, w)X_{ij}(r_1, z_1)X_{ij}(r_2, z_2)$
are expressed in terms of their product expansions similarly by
using Lemma 2.29.
Note that for all $i, j, k, l\in \Bbb Z$
$$\align & : X_{ij}(r_1, z_1)X_{ij}(r_2, z_2)X_{kl}(r_3, w):=
 :X_{kl}(r_3, w) X_{ij}(r_1, z_1)X_{ij}(r_2, z_2) :\\
&\qquad= (-1)^{(\ep_i -\ep_j, \ep_k -\ep_l)}: X_{ij}(r_1, z_1)X_{kl}(r_3, w)X_{ij}(r_2, z_2) : .\\
\endalign $$

Case 1. For $0\leq i \leq N-2$, using (3.10) and its associates
as well as Lemma 3.9 we have
$$\align & 
\{E_i(z_1)E_i(z_2) E_{i+1}(w)
-(q+ q^{-1})E_{i}(z_1)E_{i+1}(w)E_{i}(z_2)  \\
 &\qquad +  E_{i+1}(w)E_{i}(z_1)E_{i}(z_2)\}+
\{z_1 \leftrightarrow z_2 \} \\
 =& :X_{i, i+1}(0, p^{\frac{1}{2}}d^iz_1)X_{i, i+1}(0, p^{\frac{1}{2}}d^iz_2) 
 X_{i+1, i+2}(0, p^{\frac{1}{2}}d^{i+1}w):S(z_1, z_2, dw) \\
 =& 0
 \endalign $$
\medskip

Case 2. It follows from (3.10) that
$$\align &X_{12}(0, p^{\frac{1}{2}}dz_1)X_{12}(0, p^{\frac{1}{2}}dz_2) X_{01}(1, p^{-\frac{1}{2}}w)\\
=&: X_{12}(0, p^{\frac{1}{2}}dz_1)X_{12}(0, p^{\frac{1}{2}}dz_2) X_{01}(1, p^{-\frac{1}{2}}w):\\
&\quad\cdot\frac{p^{-1}w}{(z_1z_2)^{1/2}}
\frac{(z_1-z_2)(z_1-q^{-2}z_2)}{(z_1-q^{-1}d^{-1}w)(z_2-q^{-1}d^{-1}w)}
\endalign $$
Then we have 
$$\align 
& \{E_1(z_1)E_1(z_2) E_{0}(w)
-(q+ q^{-1})E_{1}(z_1)E_{0}(w)E_{1}(z_2)
  + E_{0}(w)E_{1}(z_1)E_{1}(z_2) \}\\
 &+\{ z_1 \leftrightarrow z_2 \} \\
 =& :X_{12}(0, p^{\frac{1}{2}}dz_1)X_{12}(0, p^{\frac{1}{2}}dz_2) X_{01}(1, p^{-\frac{1}{2}}w): S(z_1, z_2, d^{-1}w)p^{-1} \\
 =& 0
 \endalign $$

Case 3. In the case of $(i, i, j)=(0, 0, 1)$ we have as above

$$\align & \{E_0(z_1)E_0(z_2) E_{1}(w)
-(q+ q^{-1})E_{0}(z_1)E_{1}(w)E_{0}(z_2)
  + E_{1}(w)E_{0}(z_1)E_{0}(z_2)\} \\
 &+\{ z_1 \leftrightarrow z_2 \} \\
 =& :X_{01}(1, p^{-\frac{1}{2}}z_1)X_{01}(1, p^{-\frac{1}{2}}z_2) X_{12}(0, p^{\frac{1}{2}}dw):S(z_1, z_2, dw) \\
 =& 0
 \endalign $$

Case 4. As above we have 
$$\align 
& \{E_{N-1}(z_1)E_{N-1}(z_2) E_{0}(w)
-(q+ q^{-1})E_{N-1}(z_1)E_{0}(w)E_{N-1}(z_2)\\
 &\qquad + E_{0}(w)E_{N-1}(z_1)E_{N-1}(z_2)\} +\{ z_1 \leftrightarrow z_2 \} \\
=& :X_{N-1, N}(0, p^{\frac{1}{2}}d^{N-1}z_1)X_{N-1, N}(0, p^{\frac{1}{2}}d^{N-1}z_2) 
X_{01}(1, p^{-\frac{1}{2}}w):\\
 &\cdot S(z_1, z_2, dw) \\
 =& 0,
 \endalign $$
where we have used $p^{-1}d^{1-N}=d$.
 
Case 5. Finally we have 
$$\align 
& \{E_{0}(z_1)E_{0}(z_2) E_{N-1}(w)
-(q+ q^{-1})E_{0}(z_1)E_{N-1}(w)E_{0}(z_2)\\
 & + E_{N-1}(w)E_{0}(z_1)E_{0}(z_2)\}+\{ z_1 \leftrightarrow z_2 \} \\
=& :X_{01}(1, p^{-\frac{1}{2}}z_1)X_{01}(1, p^{-\frac{1}{2}}z_2) X_{N-1, N}(0, p^{\frac{1}{2}}d^{N-1}w):\\
 &\cdot S(z_1, z_2, pd^{N-1}w)p \\
 =& 0
 \endalign $$
 
The quantum Serre relations involving $F_i$'s are shown in a similar
way.     
Let $\alpha_i$ and $\alpha_j$ be adjacent simple roots in the affine
Dynkin diagram of type $A_{{}_{N-1}}^{(1)}$. Corresponding to the five
cases we have considered above, we have
$$\align 
& \{F_i(z_1)F_i(z_2) F_{j}(w)
-(q+ q^{-1})F_{i}(z_1)F_{j}(w)F_{i}(z_2)
  + F_{j}(w)F_{i}(z_1)F_{i}(z_2)\} \\
&\qquad\qquad +\{ z_1 \leftrightarrow z_2 \} \\
=&:F_i(z_1)F_i(z_2) F_{j}(w): 
\cases S(z_1, z_2, dw) & j=i+1\leq N-1\\
S(z_1, z_2, d^{-1}w) & i=1, \quad j=0\\
S(z_1, z_2, dw)p^{-1} & i=0, \quad j=1\\
S(z_1, z_2, d^2w)p & i=N-1, \quad  j=0\\
S(z_1, z_2, d^{-1}w)p & i=0, \quad j=N-1
\endcases\\
=& 0
\endalign
$$

 Therefore the proof is completed.\qed
\enddemo

\proclaim{Theorem 3.11} If $pq^{\pm N}$ is not a root of unity,
then 
$V_Q$ is an
irreducible  $\Cal A$ module. Therefore $V_Q$ is an irreducible $\tor$ module.
\endproclaim

\demo{Proof}   Let $U$ be a nonzero submodule of 
$$ V_Q =S(\Cal{H}^-)\otimes_\bc \bc[Q].\tag 3.12$$

Consider the subalgebra $\Cal M$ generated by the coefficient operators of $k_i^{\pm }(z)$.
Since
$$\exp(\sum_{n=1}^\infty x_n z^n ) =\sum_{n=0}^\infty y_n z^n \tag 3.13$$
where 
$$\align & y_0 =1, y_1 = x_1, y_2 = x_2 + \frac{x_1^2}{2}, y_3 = x_3 + x_2x_1 + \frac{x_1^3}{6},\\
& y_4 = x_4 + x_3 x_1 + \frac{x_2x_1^2}{2} + \frac{x_2^2}{2} + \frac{x_1^4}{24}, \cdots, \endalign $$
one can easily see that the algebra $\Cal M$ is the same as the algebra 
generated by $K_i^{\pm}=q^{\pm(\ep_i -\ep_{i+1})} $ and 
$$\align &H_{i, n}= q^{\frac{1}{2}|n|}\ep_i(n)-q^{-\frac{1}{2}|n|}\ep_j(n), 
\quad 1\leq i \leq N-1, n\in\Bbb{Z}\setminus\{0\}, \tag 3.14 \\
&H_{0, n} = q^{\frac{1}{2}|n|}\ep_N(n)-p^{-n}q^{-\frac{1}{2}|n|} \ep_1(n), \quad n\in\Bbb{Z}\setminus\{0\}.
\tag 3.15
\endalign$$
from (2.11)-(2.12) we get
$$\ep_{1N}(n) = q^{\frac{N-1}{2}|n|}\ep_1(n)-q^{\frac{1-N}{2}|n|}\ep_N(n) \in \Cal{M} \tag 3.16$$
It follows from (3.15), (3.16) and our assumption on $p$
that 
$$\ep_1(n), \ep_N(n) \in \Cal{M}, \text{ for } n\in\Bbb{Z}\setminus\{0\}.$$
Therefore $\Cal M$ contains the Heisenberg algebra $\Cal H$.
 Lemma 9.13 in \cite{K} (or Theorem 1.7.3 in \cite{FLM})
 implies that $U$ is completely reducible as  $\Cal H$-module and so
$$U= S(\wh{\Cal H}^-)\otimes \Omega$$
for some subspaces $\Omega$ of $\bc[Q]$. Suppose that 
$f=\sum_{i=1}^m s_i e^{\gamma_i} \in \Omega$, where 
 $s_i\in\bc, s_i\neq 0$, $\gamma_i\in Q$, for $1\leq i\leq m$, 
and $\gamma_i \neq \gamma_j$ if $i\neq j$.

We claim that $e^{\gamma_k}\in\Omega $ for some $k$. 

Since $K_i^{\pm} = q^{\pm(\ep_i - \ep_{i+1})}$ lie in $\Cal M$, clearly, $\Cal M$ contains $q^Q$.
Pick $\alpha\in Q$ such that $(\alpha, \gamma_{m-1}-\gamma_m)\neq 0$. Then
$$q^\alpha f - q^{(\alpha, \gamma_m)}f 
=\sum_{i=1}^{m-1}s_i( q^{(\alpha, \gamma_i)}- q^{(\alpha, \gamma_m)})e^{\gamma_i}=
\sum_{i=1}^{m-1}s_i^\prime e^{\gamma_i} \in \Omega,$$
where $s_i^\prime = s_i(q^{(\alpha, \gamma_i)} -q^{(\alpha, \gamma_m)})$, $1\leq i\leq m-1$, and 
$q^\alpha e^\gamma = q^{(\alpha, \gamma)}e^\gamma$. Since
$s_{m-1}^\prime\neq 0$ (thanks to the generic $q$), we may continue this process and get
 some $e^{\gamma_k}\in\Omega$.

Also, from (2.21), we have
$$\align &e_{\ep_i-\ep_{i+1}}\\
= &\exp(\sum_{n\in -\bz_+}\frac{\ep_i(n)-q^{n}\ep_{i+1}(n)}{n}z^{-n})
X_{i, i+1}(0, z)\\
&\cdot \exp(\sum_{n\in\bz_+}\frac{\ep_i(n)-q^{-n}\ep_{i+1}(n)}{n}z^{-n})
z^{\ep_{i+1}-\ep_i -1}\endalign$$
for $1\leq i \leq N-1$. Similarly, for $e_{\ep_{i+1}-\ep_i}$. It then follows that 
$e^{\gamma_k +Q}\subseteq \Omega$ and so $\Omega=\bc[Q]$. \qed\enddemo

\proclaim{Corollary 3.17} $V_Q$ is a level-1 $\dot{U}_v =\sln$ module
and also a level-0 $\dot{U}_h =\sln$ module.
\endproclaim

Let $\Cal B$ be the unital associative algebra  generated by the coefficients of $E_i(z),
F_i(z)$ for $1\leq i \leq N-1$ (we may assume that $d=p=1$) 
and $\ep_i(n), q^{\pm \ep_i}, 1\leq i \leq N, 
n\in\Bbb{Z}\setminus\{0\}$. Note that $\Cal M$ is a subalgebra of $\Cal B$. Then it follows 
from the proof of Theorem 3.16, we have

\proclaim{Proposition 3.18} $V_Q$ is an irreducible  $\Cal B$ module.
\endproclaim

Our next business is to show that $\Cal B$ is a homomorphic image of 
$\gln$. Thus $V_Q$ is an irreducible 
$\gln$-module. 

\medskip

\subhead \S 4.  The quantum affine algebras 
$\gln$ and $\sln$
\endsubhead

\medskip

In this section we will give a new realization for the quantum affine
algebra $\gln$.

We define the quantum affine algebra $\gln$ to be the 
unital associative algebra generated by $\ep_{im}, q^{c/2}, k_{i0}^{\pm 1}$, $x^{\pm}_{jn}$ $(i=1, \cdots, N; j=1, \cdots, N-1; n\in\Bbb Z, m\in\Bbb Z^{\times})$ subject to 
the following relations
%that $x^{\pm}_i(z)=\sum_n x^{\pm}_{in}x^{-n}$,
that $q^{c/2}$ is central and
$$\align
[\ep_{im}, \ep_{jn}]&=\frac{[m]}{m}[mc]\delta_{ij}\delta_{m, -n}, \tag4.1\\
[\ep_{im}, x^{\pm}_{jn}]&=0, \qquad j\neq i, i-1, \tag4.2\\
[\ep_{im}, x^{\pm}_{in}]&=\pm \frac{[m]}{m}q^{\mp |m|c/2+|m|/2}x^{\pm}_{i,m+n}, \tag4.3\\
[\ep_{im}, x^{\pm}_{i-1,n}]&=\mp \frac{[m]}{m}q^{\mp |m|c/2-|m|/2}x^{\pm}_{i-1,m+n}, \tag4.4\\
x^{\pm}_{i,m+1}x^{\pm}_{jn}&-q^{\pm a_{ij}}x^{\pm}_{jn}x^{\pm}_{i, m+1}
=q^{\pm a_{ij}}x^{\pm}_{im}x^{\pm}_{j,n+1}-x^{\pm}_{j,n+1}x^{\pm}_{im},
\tag4.5\\
[x^{\pm}_{im},x^{\pm}_{jn}]&=0, \qquad |j-i|>1, \tag4.6\\
[x^+_{im}, x^-_{jn}]&=\frac{\delta_{ij}}{q-q^{-1}}
(k^+_{i, m+n}q^{(m-n)c/2}-
k^-_{i, m+n}q^{(n-m)c/2}), \tag4.7\\
x^{\pm}_{im_1}x^{\pm}_{im_2}x^{\pm}_{i\pm 1, n}&-(q+q^{-1})x^{\pm}_{im_1}x^{\pm}_{i\pm 1, n}x^{\pm}_{im_2} \\
\qquad & +x^{\pm}_{i\pm 1, n}x^{\pm}_{im_1}x^{\pm}_{im_2}+\{m_1\leftrightarrow m_2\}=0, \tag4.8
\endalign
$$
where $(a_{ij})$ is the Cartan
matrix of type $A$ of rank $N-1$, $[m]=\frac{q^m-q^{-m}}{q-q^{-1}}$,  
$[mc]=\frac{q^{mc}-q^{-mc}}{q-q^{-1}}$, and 
$$
k_i^{\pm}(z)=\sum_{n=0}^\infty k^{\pm}_{in}z^{-n}=k_{i0}^{\pm 1}exp(\pm (q-q^{-1})\sum_{n>0}
(q^{n/2}{\ep_{in}}-q^{-n/2}\ep_{i+1, n})z^{\pm n}). \tag4.9
$$

Let $\ep_i(m)=\frac{m}{[m]}\ep_{im}$ and $c=1$, then we have
$$
[\ep_i(m), \ep_j(n)]=m\delta_{ij}\delta_{m, -n}I,
$$
Thus the Fock module $V_Q$ constructed in Sect. 3 is a
vertex representation of the quantum affine algebra $\gln$
at level one. 

The following easily checked result 
shows that $\gln$ contains the quantum affine
algebra $\sln$ canonically.
\proclaim{Lemma 4.10} Let $h_{im}=q^{|m|/2}\ep_{im}-q^{-|m|/2}\ep_{i+1,m}$, then the associative subalgebra
generated by $h_{im}, q^{c/2}, k_{i0}^{\pm 1}$, $x^{\pm}_{jn}$ 
$(i,j=1, \cdots, N-1; n\in\Bbb Z, m\in\Bbb Z^{\times})$
is isomorphic to the quantum affine algebra $\sln$.
The commutation relations are
as follows: Eqs. (4.5)-(4.9) and 
$$
\align
[h_{im}, h_{jn}]&=\frac{[a_{ij}m]}{m}[mc]\delta_{m, -n}, \tag4.11\\
[h_{im}, x^{\pm}_{jn}]&=\pm \frac{[a_{ij}m]}{m}q^{\mp |m|c/2}x^{\pm}_{j,m+n}, \tag4.12
\endalign
$$
\endproclaim

\proclaim{Remark 4.13}
In terms of generating functions $x^{\pm}_i(z)=\sum_n 
x^{\pm}_{in}z^{-n}$,
$k^{\pm}_i(z)$, the commutation relations of the quantum affine
algebra $\sln$ are given in Eqs. (1.2)-(1.13)
with $e_i(z)=x_i^+(z), f_i(z)=x^-(z)$ and $d=1$.
\endproclaim

We now introduce another set of basis in the Heisenberg subalgebra
to establish isomorphism between the quantum affine algebra $\gln$
defined above
and the Ding-Frenkel algebra $\gln$ \cite{DF}.

For $m\neq 0$ we define
$$\align
&a_{Nm}\\
&=q^{Nm+|m|/2}\left(\frac{|m|}m
\frac{\ep_{1m}+q^{|m|}\ep_{2m}+\cdots+q^{(N-1)|m|}
\ep_{Nm}}{(1+q^{2|m|}+\cdots+q^{2(N-1)|m|})^{1/2}}+\ep_{Nm}\right). \tag4.14
\endalign
$$
For each $i=1, \cdots, N-1$ let 
$$
a_{im}=\sum_{j=i}^{N-1}q^{im}(q^{|m|/2}\ep_{im}-q^{-|m|/2}\ep_{i+1, m})
+a_{Nm}. \tag4.15
$$

\proclaim{Theorem 4.16} The quantum affine algebra $\gln$
is isomorphic to the associative algebra generated by the generators  $ x^{\pm}_{in}$,
$a_{im}, q^{c/2}, k_{i0}^{\pm 1}$, $x^{\pm}_{jn}$  $(i=1, \cdots, N; j=1, \cdots, N-1; n\in\Bbb Z, m\in\Bbb Z^{\times})$ with the commutation 
relations (4.5-4.8) and: 
$$
\align
[a_{im}, a_{in}]&=0, \tag 4.17\\
[a_{im}, a_{jn}]&=-\frac{[m]}m[mc]q^{-m}\delta_{m, -n},
\qquad\qquad i>j\tag 4.18\\
[a_{im}, x^{\pm}_{jn}]&=0,\qquad\qquad j\neq i, i+1,\tag 4.19\\
[a_{im}, x^{\pm}_{in}]&=\pm\frac{[m]}mq^{(i-1)m\mp|m|c/2} x^{\pm}_{i,m+n},\tag4.20\\
[a_{i+1,m}, x^{\pm}_{in}]&=\mp\frac{[m]}mq^{(i+1)m\mp|m|c/2} x^{\pm}_{i,m+n},\tag4.21
\endalign
$$
\endproclaim
\demo{Proof} Let $\oh_{im}=q^{im}h_{im}=q^{im}(q^{|m|/2}\ep_{im}-q^{-|m|/2}\ep_{i+1,m})$. First we compute that
$$\align
[&a_{Nm}, a_{Nn}]=\frac{[m][mc]}m \delta_{m,-n}\\ 
&\cdot\big(-q^{|m|}+(\frac{|m|}m+\frac{|n|}n)\frac{1}{1+q^{2|m|}+
\cdots +q^{2(N-1)|m|}}+q^{|m|})=0,\tag4.22
\endalign
$$
and we can also easily check that
$$
[\oh_{im}, a_{Nn}]=-\delta_{i, N-1}\frac{q^{-m}}m[m][mc]\delta_{m, -n}.
\tag4.23
$$

It follows from Lemma 4.10 and Eq.(4.23) that
$$\align
[a_{im}, \oh_{jn}]&=[\oh_{im}+\cdots+\oh_{N-1,m}+a_{Nm}, \oh_{jn}]\\
&=\cases 0, & i\neq j, j+1\\
 -\frac{q^m}m[m][mc]\delta_{m. -n}, & i=j+1\\
\frac{q^{-m}}m[m][mc]\delta_{m. -n}, & i=j \endcases . \tag4.24
\endalign
$$
We claim that the new set of generators of the Heisenberg subalgebra
satisfy the following relations:
$$
[a_{im}, a_{jn}]=\cases 0, & i=j\\
\frac {-q^m}m[m][mc]\delta_{m, -n}, &\text{for } i>j. \endcases
\tag4.25
$$
In fact for $i>j$ it follows from Eqs. (4.22-23) that
$$\align
[a_{im}, a_{jn}]&=[\oh_{im}+a_{i+1,m}, a_{jn}]\\
&=[a_{i+1, m}, a_{jn}]=\cdots=[a_{Nm}, a_{jn}]\\
&=[a_{Nm}, \oh_{jn}+\oh_{j+1,n}+\cdots+\oh_{N-1,n}+a_{Nn}]\\
&=[a_{Nm}, \oh_{N-1,n}]=\frac {-q^m}m[m][mc]\delta_{m, -n}.
\endalign
$$
Similarly using Eq. (4.23) we obtain that 
$$
[a_{im}, a_{in}]=[a_{i+1,m}, a_{i+1,n}]=\cdots=[a_{Nm}, a_{Nn}]=0.
$$
Now let's look at the commutation relations
$[a_{im}, x^{\pm}_{jn}]$. For simplicity we write
$$a_{Nm}=c_{1m}\ep_{1m}+c_{2m}\ep_{2m}+\cdots+
c_{Nm}\ep_{Nm}.
$$
Then it follows that  $c_{i+1,m}=q^{|m|}c_{im}$ for $i=1, \ldots, N-2$ and $c_{Nm}$=$q^{|m|/2}(c_{N-1, m}q^{|m|/2}+q^{Nm})$. Using the latter equation and (4.2-4.3) we have that
$$\align
[a_{Nm}, x^{\pm}_{N-1, n}]&=\pm\frac{[m]}mq^{\mp|m|c/2}x^{\pm}_{N-1, n}
(c_{N-1, m}q^{|m|/2}-c_{Nm}q^{-|m|/2})\\
&=\mp\frac{[m]}mq^{\mp|m|c/2+Nm}x^{\pm}_{N-1, n}.
\endalign
$$
Similarly one proves that $[a_{Nm}, x^{\pm}_{in}]=0$ for
$1\leq i\leq N-2$ by using $c_{i+1,m}=q^{|m|}c_{im}$. 
The relations (4.18-19) are then immediate consequences of $a_{im}=a_{i+1, m}+q^{im}h_{im}$ and Lemma 4.10.

The other relations hold automatically, and  the theorem is proved. \qed\enddemo

The following theorem is obtained by the standard
argument of generating functions. 

\proclaim{Theorem 4.26} The new set of generators 
$a_{im}, q^{c/2}, k_{i0}^{\pm 1}$, $x^{\pm}_{jn}$ of $\gln$ 
$(i=1, \cdots, N; $ 
$j=1, \cdots, N-1; n\in\Bbb Z, m\in\Bbb Z^{\times})$ satisfy the following
commutation relations in terms of generating functions $v_i^{\pm}(z)$, $x_i^{\pm}(z)$ where
$$v_i^{\pm}(z)=v_{i0}^{\mp 1}exp(\mp (q-q^{-1})\sum_{n>0}{a_{i, \mp n}}z^{\pm n}),
$$
and the relations are
$$
\align
v_i^{\pm}(z)v_j^{\pm}(w)&=v_j^{\pm}(w)v_i^{\pm}(z), \qquad 
v^{\pm}_{j0}v_{j0}^{\mp}=1, \tag 4.27\\
v_i^{+}(z)v_i^{-}(w)&=v_i^{-}(w)v_i^{+}(z), \tag 4.28\\
\frac{z-wq^{\pm c}}{z-wq^{\pm c-2}}v_i^{\mp}(z)v_j^{\pm}(w)
&=v_j^{\pm}(w)v_i^{\mp}(z)\frac{z-wq^{\mp c}}{z-wq^{\mp c-2}}, \quad j>i, \tag4.29\\
v_i^{\pm}(z)^{-1}x_j^{\pm}(w)v_i^{\pm}(z)&=x_j^{\pm}(w), \quad i-j\leq -1, \tag4.30\\
v_i^{\pm}(z)^{-1}x_j^{\mp}(w)v_i^{\pm}(z)&=x_j^{\mp}(w), \quad i-j\geq 2, \tag4.31\\
v_i^{\pm}(z)^{-1}x_i^{\vep}(w)v_i^{\pm}(z)&=
\frac{zq^{\mp \vep c/2+1}-wq^{-i-1}}{zq^{\mp \vep c/2}-wq^{-i}}x_i^{\vep}(w), \quad\vep=+ \ \text{or} \ - \tag4.32\\
v_{i+1}^{\pm}(z)^{-1}x_i^{\vep}(w)v_{i+1}^{\pm}(z)&=
\frac{zq^{\mp \vep c/2-1}-wq^{-i+1}}{zq^{\mp \vep c/2}-wq^{-i}}x_i^{\vep}(w), 
\quad \vep=+ \ \text{or} \ - \tag4.33\\
(z-wq^{\mp a_{ij}})x_i^{\pm}(z)x_j^{\pm}(w)&=x_j^{\pm}(w)x_i^{\pm}(z)
(zq^{\mp a_{ij}}-w), \tag 4.34
\endalign
$$
$$\align
[x_i^{\pm}(z), x_j^{\pm}(w)]&=0, \quad |i-j|>1, \tag4.35\\
[x_i^{+}(z), x_j^{-}(w)]&=\frac{\delta_{ij}}{q-q^{-1}}\big\{\delta(\frac z{wq^c})v_{i+1}^-(wq^{c/2})v_i^-(wq^{c/2})^{-1}\\
&\qquad \qquad -\delta(\frac {zq^c}w)v_{i+1}^+(zq^{c/2})v_i^+(zq^{c/2})^{-1}\big\}, \tag4.36\\
\big\{x_i^{\pm}(z_1)x_i^{\pm}(z_2)x_j^{\pm}&(w)-(q+q^{-1})x_i^{\pm}(z_1)x_j^{\pm}(w)x_i^{\pm}(z_2) \\
+x_j^{\pm}(w)&x_i^{\pm}(z_1)x_i^{\pm}(z_2)\big\} +\{z_1 \leftrightarrow z_2\}=0, \quad i-j=\pm 1, \tag4.37
\endalign
$$
where the rational functions represent the corresponding Taylor expansions as usual.
\endproclaim

\proclaim{Remark 4.38} One can also express the generator $\ep_{im}$ in
terms of $a_{im}$ easily from Eqs. (4.14-4.15) to give another proof
of the theorem.
If we let
$$
X_i^{\pm}(z)=x_i^{\pm}(q^{-i}z),
$$
then the generators defined by $X^{\pm}_i(z), v_j^{\pm}(z), q^{c/2}$
are exactly the generators in Ding-Frenkel's definition of $\gln$
(with correction of some typos in \cite{DF}).  
Our Heisenberg algebra is orthogonal and a direct deformation
of the classical affine general linear algebra.

Under the presentation of Theorem 4.26 the subalgebra $\sln$
 is generated by $q^c$ and 
$$x_i^{\pm}(z), k^+_i(z)=v_{i+1}^-(zq^{i})v_i^-(zq^{i})^{-1},
 k^-_i(z)=v_{i+1}^+(zq^{i})v_i^+(zq^{i})^{-1}.$$

\medskip
\endproclaim

In order to let $V_Q$ have a weight space decomposition, we need to add the operator $q^D$ to 
$\gln$ and $\Cal B$.
As $\gln$-module, we have $V_Q= \oplus V^\mu$, where 
$$V^\mu = \{v\in V_Q: q^h.v = q^{\mu(h)}v, \text{ for } h\in \frak{h}\}.$$
It is easy to see that
$V^\mu = V_\mu$, where $V_\mu = \{ v\in V_Q : h.v = \mu(h) v, \text{ for } h\in \frak{h}\}$.

\proclaim{Proposition 4.39} $V_Q$ is an irreducible $\gln$  module 
with the weight space decomposition $V_Q= \sum_{\mu\in P}\oplus V_\mu$. Moreover, $V_Q$
is complete reducible as $\sln$-module.
\endproclaim

\enddocument